\pdfoutput=1

\documentclass[12pt]{article}
  
\usepackage{amssymb,amsmath,amsthm,geometry}
\usepackage[square,numbers]{natbib}
\usepackage{IEEEtrantools}
\usepackage{mathtools}
\usepackage{tikz}
\usepackage{subcaption}
\usepackage{mathabx}
\usepackage{hyperref}
\usetikzlibrary{calc}
\usepackage{graphicx}
\usepackage{multirow}
\usepackage{enumerate}
\usepackage[title]{appendix}

\usepackage{graphics}

\newtheorem{innercustomgeneric}{\customgenericname}
\providecommand{\customgenericname}{}
\newcommand{\newcustomtheorem}[2]{%
  \newenvironment{#1}[1]
  {%
   \renewcommand\customgenericname{#2}%
   \renewcommand\theinnercustomgeneric{##1}%
   \innercustomgeneric
  }
  {\endinnercustomgeneric}
}

\newcustomtheorem{customthm}{Theorem}
\newcustomtheorem{customlemma}{Lemma}
\newcustomtheorem{customcor}{Corollary}

\DeclareMathOperator{\argmin}{argmin}

\hypersetup{
    colorlinks=true,
    linkcolor=blue,
    filecolor = blue
, urlcolor = blue,
citecolor = blue
}

\title{Equivariant Estimation of Fr\'echet Means} 
\author{Andrew McCormack and Peter Hoff \\
Department of Statistical Science \\
Duke University} 
\date{\today}

\begin{document}
\maketitle

\begin{abstract}
The Fr\'echet mean generalizes the concept of a mean to a metric space setting. In this work we consider equivariant estimation of Fr\'echet means for parametric models on metric spaces that are Riemannian manifolds. The geometry and symmetry of such a space is encoded by its isometry group. Estimators that are equivariant under the isometry group take into account the symmetry of the metric space. For some models there exists an optimal equivariant estimator, which necessarily will perform as well or better than other common equivariant estimators, such as the maximum likelihood estimator or the sample Fr\'echet mean. We derive the general form of this minimum risk equivariant estimator and in a few cases provide explicit expressions for it. In other models the isometry group is not large enough relative to the parametric family of distributions for there to exist a minimum risk equivariant estimator. In such cases, we introduce an adaptive equivariant estimator that uses the data to select a submodel for which there is an MRE. Simulations results show that the adaptive equivariant estimator performs favorably relative to alternative estimators.

\end{abstract}

\section{Introduction}
Data analysis settings where observations do not take values in a vector space present unique challenges. One such setting is metric space $(\mathcal{X},d)$ valued data where there is not necessarily a notion of addition and scalar multiplication of points in $\mathcal{X}$, but there is a notion of distance between points. Classical examples of metric space-valued data are directional data such as sphere-valued data, orthonormal frame data and subspace data \cite{mardia2009directional,chikuse2012statistics}. Detailed expositions of metric space structures on the Stiefel manifold $\text{V}_{k}(\mathbb{R}^n)$ of orthonormal $k$-frames in $\mathbb{R}^n$ and the real Grassmannian manifold $\text{Gr}_{k}(\mathbb{R}^n)$ of $k$-dimensional subspaces in $\mathbb{R}^n$ can be found in \cite{edelman1998geometry, lim1807grassmannian}. Other notable examples of metric space-valued data include positive definite covariance matrices \cite{lin2019riemannian}, shape space modelling on the complex Grassmannian \cite{kendall1984shape,kent1994complex}, and hierarchical structures that can be represented in hyperbolic spaces \cite{nagano2019wrapped}. 

 If $\mathcal{X}$ is not a vector space then the arithmetic mean is not available as a description of location. However, the distance function describes the relative locations of points in a metric space and allows the notion of a mean to be generalized. For $k > 0$, the $k$-\textit{Fr\'echet mean} \cite{frechet1948elements} of the law $P$ of a  metric space-valued random object $X \sim P$ taking values in $(\mathcal{X},d)$, is defined as
\begin{align}
    E_kX = E_kP \coloneqq \underset{x \in \mathcal{X}}{\text{argmin}} \; E\big(d(X,x)^k \big).
    \label{frecmeandefinition}
\end{align}
In words, a $k$-Fr\'echet mean is the collection of points in $\mathcal{X}$ that are on average the closest to $X$ with respect to the $k$th power of the distance function. When $\mathcal{X} = \mathbb{R}$ under the Euclidean metric, $E_1X$ and $E_2X$ correspond to the usual median and mean respectively.  Just like medians in $\mathbb{R}$, $k$-Fr\'echet means are set valued. The $2$-Fr\'echet mean is the primary estimand of interest in this article and will be denoted by $EX$ or $EP$ and will be referred to as the Fr\'echet mean.

 The most basic nonparametric estimator of $EP$ given realizations $x_1,\ldots,x_n$ of i.i.d.\@ random objects $X_1,\ldots,X_n$ distributed according to $P$, is the sample Fr\'echet mean $\Bar{X}$, defined by 
\begin{align}
    \Bar{X} \coloneqq E_2\bigg( \frac{1}{n}\sum_{i = 1}^n \delta_{x_i}\bigg) = \underset{x \in \mathcal{X}}{\text{argmin}} \; \frac{1}{n} \sum_{i = 1}^n d(x_i,x)^2, \label{sampfrecmn}
\end{align}
where $\delta_{X_i}$ is a Dirac measure. Being an $M$-estimator, the convergence properties of $\Bar{X}$ to $EP$ are non-trivial and are of substantial interest, especially in relation to how the geometry of $\mathcal{X}$ impacts rates of convergence \cite{schotz2019convergence,eltzner2019smeary,pennec2019curvature}. However, in a parametric setting the sample Fr\'echet mean may not be the most efficient estimator as it does not utilize any information about the parametric family under consideration. It is the goal of this article to provide alternative estimators of $EP$ in parametric settings using ideas from equivariant estimation theory.  

Parametric models for metric space-valued data are typically tailored to the metric space $\mathcal{X}$. Many of the classical models for directional data such as the Langevin, Bingham and von Mises-Fisher distributions are exponential families. A method for constructing exponential families on a homogeneous manifold is provided in \cite{tojo2019method,cohen2015harmonic} building on work on exponential transformation models \cite{barndorff1982exponential}. Estimation of Fr\'echet means may also be of interest for less exotic spaces if they are endowed with a non-standard metric. An example of this is provided in Section \ref{SecNontrans} where the Fr\'echet mean of a Wishart-distributed matrix is estimated using the log-Euclidean metric as a loss function \cite{arsigny2006log}. General parametric models that apply to arbitrary metric spaces are less common with a notable exception being the Riemannian Gaussian distribution \cite{pennec2006intrinsic}.

Parameter estimation in models of directional data have largely focused on large sample asymptotics, maximum likelihood methods and Bayesian methods \cite{chikuse2012statistics,  pal2020conjugate}. The problem of specifically estimating a Fr\'echet mean has not been treated  extensively from a decision-theoretic perspective. Recent work in this area has considered the properties of shrinkage estimators for Fr\'echet means \cite{yang2019shrinkage,yang2020shrinkage,mccormack2020stein}. These works highlight that estimators such as the MLE or sample Fr\'echet mean can be inadmissible. As an alternative to improving upon $\Bar{X}$ or the MLE via shrinkage, one could consider finding the optimal equivariant estimator within a class of equivariant estimators. Specifically, every metric space inherits a group that preserves its metric structure, the isometry group of distance preserving bijections of $\mathcal{X}$. These isometries act on $\mathcal{X}$ and give the Fr\'echet mean estimation problem the structure of an invariant decision problem. If it exists, the optimal equivariant estimator for this problem will perform at least as well as the sample Fr\'echet mean or MLE, since these are both isometrically equivariant. The minimum risk equivariant estimator of a Fr\'echet mean can be seen as a natural generalization of the Pitman estimator for location families on the real line \cite{pitman1939estimation}. It is also a proper Bayes estimator if the isometry group is compact, in which case it must also be admissible.

An outline of this article is as follows: Section \ref{Secprelim} reviews relevant group theoretic concepts. Section \ref{SecEquivariantEst} introduces the equivariant estimation problem and provides a form for the optimal equivariant estimator under a transitive isometry group action. The optimal equivariant estimator can be viewed as a Bayes estimator under a prior induced by the right Haar measure on the isometry group. It is shown that under certain conditions, the optimal equivariant estimator can alternatively be characterized as a Bayes estimator where a uniform prior is place on the Fr\'echet mean. Explicit expressions for the optimal equivariant estimator are derived for generalizations of the von Mises-Fisher distributions on the sphere, hyperbolic space and Stiefel manifold. In these cases the MLE is equal to the MRE. Section \ref{SecNontrans} explores the more common scenario where the isometry group does not act transitively. In this case we propose an estimator that adaptively selects a submodel where the isometry group is transitive from which an equivariant estimator is constructed. Simulation studies on the space of positive definite matrices and the $p$-torus illustrate the efficacy of this adaptive equivariant estimator.

\section{Mathematical Preliminaries}
\label{Secprelim}

Associated with a metric space $(\mathcal{X},d)$ is the set of bijections from $\mathcal{X}$ to itself that preserve distances, namely the isometry group
\begin{align}
\label{isometrydefinition}
    \text{Iso}(\mathcal{X}) \coloneqq \{g: g \;\; \text{surjective},\;  d(g(x),g(y)) = d(x,y), \; \forall x,y \in \mathcal{X}\}.
\end{align}
The isometry group is a group under function composition and acts on $\mathcal{X}$ via the evaluation map $gx \coloneqq g(x)$. To ease notation $G$ will be used interchangeably with $\text{Iso}(\mathcal{X})$ throughout this article. The \textit{isotropy} group at $x_0$ of a group $G$ acting on a space $\mathcal{X}$ is the subgroup of $G$ defined by 
\begin{align}
\label{isotropydefinition}
    G_{x_0} \coloneqq \{g: gx_0 = x_0\}.
\end{align} Any action of a group $G$ on $\mathcal{X}$ partitions $\mathcal{X}$ into equivalence classes determined by the equivalence relation $x \sim y$ if there exists a $g \in G$ with $gx = y$. Each such equivalence class is called an \textit{orbit} and the orbit containing $x \in \mathcal{X}$ will be denoted by $[x]$, with $\mathcal{X}/G$ denoting the collection of all orbits of this $G$ action. If there is only one orbit, $G$ is said to   act \textit{transitively} on $\mathcal{X}$.  If $\text{Iso}(\mathcal{X})$ acts transitively on $\mathcal{X}$ then $\mathcal{X}$ is called \textit{homogeneous} and possesses a high degree of symmetry, as every point $x \in \mathcal{X}$ ``looks the same" as every other point $y \in \mathcal{X}$ with respect to the metric. In this work we restrict the metric spaces under consideration to be homogeneous. Moreover, we make the extra assumption that $\mathcal{X}$ is a Riemannian manifold with Riemannian distance function $d$. This assumption is not strictly necessary in what follows but it ensures that $\text{Iso}(\mathcal{X})$ is sufficiently well-behaved as a topological group.  Additional details on the basics of Riemannian manifolds and their distance functions can be found in  \cite{lee2013smooth, carmo1992riemannian}.

In $\mathbb{R}^n$ the Lebesgue measure $\lambda$ is invariant under the group operation of addition as $\lambda(A) = \lambda(A + b)$ for any Borel $A$ and $b \in \mathbb{R}$. A generalization of Lebesgue measure to a topological group $G$ with its Borel $\sigma$-algebra are the left and right invariant \textit{Haar measures}, $\lambda_L$ and $\lambda_R$ which satisfy
\begin{align*}
    \lambda_L(gA) = \lambda_L(A), \;\; \lambda_R(Ag) = \lambda_R(Ag), \;\; \forall A \in \mathcal{B}(G), \; g \in G,
\end{align*}
where $gA \coloneqq \{gh: h \in A\}$ and similarly for $Ag$. If the topology of $G$ is locally compact then $\lambda_L$ and $\lambda_R$ exist and are unique up to scaling \cite{nachbin1976haar}. The measures $\lambda_L$ and $\lambda_R$ need not be the same. However, in an abelian group like $\mathbb{R}^n$ or a compact group like the orthogonal group $\text{O}(n)$, the left and right Haar measures agree up to scaling. The modular function is a continuous homomorphism $\Delta:G \rightarrow \mathbb{R}^{\times}$ such that $\lambda_L(Ag) = \Delta(g)\lambda_L(A)$ for all $A \in \mathcal{B}(G)$ A group is said to be \textit{unimodular} if its left and right Haar measures agree or equivalently $\Delta(g) = 1, \; \forall g \in G$. As the measure $\Delta(g)^{-1} \lambda_L(dg) = \Delta(g^{-1})\lambda_L(dg)$ can be seen to be right invariant, $\lambda_L$ and $\lambda_R$ are related by $\lambda_L(dg) = \Delta(g) \lambda_R(dg)$, again up to scaling. More generally, it is of interest to consider measures on topological spaces $\mathcal{X}$ that are acted on continuously by a topological group $G$. A measure $\nu$ on $(\mathcal{X},\mathcal{B}(\mathcal{X}))$ is \textit{relatively invariant} with multiplier $\chi(g)$ if $\nu(gA) = \chi(g)\nu(A)$ for all $A \in \mathcal{B}(\mathcal{X}), \; g \in G$. The multiplier $\chi:G \rightarrow \mathbb{R}^{\times}$ is a continuous homomorphism. Further details on the interplay between groups and measures can be found in \cite{wijsman1990invariant, eaton1989group,nachbin1976haar}.  Hausdorff measures on a metric space are relatively invariant with respect to the isometry group action with multiplier $\chi = 1$. For a Riemannian manifold $\mathcal{X}$, the Hausdorff measure is the same as the Riemannian volume measure $\text{vol}(dx)$. If $\mathcal{X}$ is a manifold embedded in $\mathbb{R}^n$ then $\text{vol}(dx)$ can be thought of as the ``surface area" measure of $\mathcal{X}$.

\section{Estimation Under a Transitive Action}
\label{SecEquivariantEst}
    \subsection{Equivariant Estimation}
An equivariant estimation problem consists of a family of distributions $\mathcal{P} = \{P_\theta:\theta \in \Theta\}$ on the sample space $\mathcal{X}$, an invariant loss function $L(\delta(x),\theta): \mathcal{D} \times \Theta \rightarrow \mathbb{R}^+$ and a group $G$ that acts measurably on $\mathcal{X}$ \cite{berger2013statistical}. For simplicity, let the decision space $\mathcal{D}$ be the same as $\Theta$. The decision problem is to estimate $\Theta$, or a functional thereof, given an observation $X \sim P_\theta$. The group $G$ induces an action on the set of all probability measures on $\mathcal{X}$ given by $P \rightarrow gP$ where $(gP)(A) \coloneqq P(g^{-1}(A))$ for all measurable sets $A$. It is assumed that $\mathcal{P}$ is invariant under this $G$ action, meaning that $gP_\theta \in \mathcal{P}$ for all $\theta,g$. It is always possible to find a family of distributions that contains $\mathcal{P}$ and is invariant, namely $G\mathcal{P} \coloneqq \{gP_\theta:g \in G,\; \theta \in \Theta\}$. It is also assumed that the parameterization $\Theta$ of $\mathcal{P}$ is identifiable so that there exists a unique $g\theta \in \Theta$ with $gP_\theta = P_{g\theta}$ and hence $G$ also acts on $\Theta$. The loss is defined to be invariant if $L(g\theta,g\delta) = L(\theta,\delta)$ where $g\delta$ and $g\theta$ are the results of $G$ acting on $\Theta$.  

An equivariant estimator is a function $\delta:\mathcal{X} \rightarrow \Theta$ satisfying $\delta(gx) = g\delta(x)$. Due to both the invariance of the loss and the invariance of the family $\mathcal{P}$, the risk function of any equivariant estimator is constant on $\Theta$-orbits:
\begin{align}
   R(\theta,\delta) =   E_{\theta}\big(L(g\theta,g\delta(X))\big) =  E_{\theta}\big(L(g\theta,\delta(gX))\big) =  R(g\theta,\delta).
\end{align}
Consequently, if $G$ acts transitively on $\Theta$ the risk functions of equivariant estimators can be totally ordered since they are constant. It is then of interest to search for the minimum risk equivariant estimator (MRE). As many standard estimators are equivariant, equivariant estimation procedures can be motivated as a way to construct estimators that outperform such commonly used estimators. When finding the MRE it can be beneficial to work with as small of a group as possible that remains transitive over $\Theta$. If $H$ is a subgroup of $G$ then any $G$ equivariant estimator is $H$ equivariant. Thus the MRE under $H$ will perform at least as well as the MRE under $G$. 

\subsection{Fr\'echet Mean Estimation Problem}
Before introducing the estimation problem, we observe that the Fr\'echet mean is an equivariant function under $G = \text{Iso}(\mathcal{X})$ \cite{chakraborty2020manifoldnet}. If $L^2(\mathcal{X})$ is the collection probability measures on $X$ with $\int d(x,y)^2P(dy) < \infty$ for at least one $x \in \mathcal{X}$ then the Fr\'echet mean can be viewed as a function $E:L^2(\mathcal{X}) \rightarrow 2^\mathcal{X}$ where $2^\mathcal{X}$ is the power set of $\mathcal{X}$. The isometry group acts on $L^2(\mathcal{X})$ by $P \rightarrow gP$ where $gP(A) = P(g^{-1}(A))$ for every Borel set $A$.
 With the natural action of $G$ on $2^\mathcal{X}$ defined by $g \bigcup_i x_i \coloneqq \bigcup_i g(x_i) \in 2^\mathcal{X}$, the Fr\'echet mean is equivariant, meaning that $EgP = gEP$. This follows from the definition \eqref{frecmeandefinition} of $E$ since if $x \in EP$ then
 \begin{align*}
     & \int d(x,y)^2P(dy) \leq \int d(z,y)^2 P(dy) \;\;\; \forall z \in \mathcal{X},
\end{align*}
which implies
     \begin{align*}
     \int d(gx,y)(gP)(dy) =  \int d(gx,gy)^2 P(dy) \leq  \int d(z,y)^2 (gP)(dy) \;\;\; \forall z \in \mathcal{X},
 \end{align*}
so that $gx \in EgP$. This proves that $gEP \subset EgP$ and applying this result with $\Tilde{P} = gP$ and $\Tilde{g} = g^{-1}$ yields $g^{-1}EgP =  \Tilde{g}E\Tilde{P} \subset E\Tilde{g}\Tilde{P} = EP$, so $EgP = gEP$ as needed. The equivariance of the Fr\'echet mean implies that if $EP$ is a singleton set then $EgP$ is also a singleton set for all $g \in G$.

The estimation problem of interest in this article is to estimate the Fr\'echet mean of $P_\theta$ under the squared distance loss function $L(EP_\theta,\delta) = d(EP_\theta,\delta)^2$, given i.i.d.\@ observations $X_1,\ldots,X_n$ from $P_\theta$. The distribution $P_\theta$ is assumed to be a member of the family of distributions $\mathcal{P} = \{P_\theta:\theta \in \Theta\}$ on a homogeneous Riemannian manifold $(\mathcal{X},d)$ where $\mathcal{P}$ is invariant under the action of $G$.  Moreover, in this section it is assumed that this action is transitive over $\Theta$. The Fr\'echet mean $EP_\theta$ is assumed to be a singleton set so that $EP_\theta \in \mathcal{X} \subset 2^\mathcal{X}$ and $d(EP_\theta,\delta)$ makes sense as a function from $\mathcal{X} \times \mathcal{X} \rightarrow \mathbb{R}$.   In practice $P_\theta$ typically has a unique Fr\'echet mean. Theoretical guarantees of the uniqueness of $EP_\theta$ usually require that $P_\theta$ be supported in a small metric ball in $\mathcal{X}$ \cite{karcher1977riemannian,kendall1990probability}. Lastly, the family $\mathcal{P}$ is assumed to be dominated by the Riemannian volume measure $\text{vol}(dx)$ on $\mathcal{X}$ where $P_\theta$ has the density $p(x|\theta)$ with respect to $\text{vol}(dx)$. As $\text{vol}(dx)$ is an invariant measure under $G$, $P_{g\theta}$ has the density $p(x|g\theta)$ and no modification by a Jacobian term is required.

The above problem is an equivariant estimation problem under the componentwise action of $G$ on $\mathcal{X}^n$ defined by $g(X_1,\ldots,X_n) \coloneqq (gX_1,\ldots,gX_n)$. By the definition of an isometry, the loss is invariant since
\begin{align*}
    L(\theta,\delta) = d(EP_\theta,\delta)^2 = d(gEP_\theta,g\delta)^2 = d(EP_{g\theta},g\delta)^2 = L(g\theta,g\delta),
\end{align*}
 and the family $\mathcal{P}_n = \{P_\theta
\times \cdots \times P_\theta:\theta \in \Theta\}$ is invariant under $G$. The collection $(X_1,\ldots,X_n)$ will be denoted by $X$ or $x$ throughout and similarly $p_n(x|\theta)  \coloneqq \prod_{i = 1}^np(x_i|\theta)$ with respect to the base measure $\text{vol}(dx) \coloneqq \bigtimes_{i = 1}^n\text{vol}(dx_i)$ on $\mathcal{X}^n$.

As shown above, Fr\'echet means are equivariant. In particular, the Fr\'echet mean under the empirical distribution, which is the sample Fr\'echet mean  \eqref{sampfrecmn}, is equivariant. The MLE defined by 
\begin{align}
    \delta_{MLE}(x) \coloneqq EP_{\hat{\theta}}, \;\;\;\; \hat{\theta}(x) \coloneqq \underset{\theta \in \Theta}{\text{argmax}} \; \sum_{i = 1}^n \log(p(x_i|\theta)), 
\end{align}
 where $\hat{\theta}(x)$ is assumed to be the unique maximizer of the log-likelihood for all $X \in \mathcal{X}^n$, is also equivariant. If the MRE exists and differs from these estimators it must necessarily outperform them in terms of squared distance loss.    

\subsection{Finding the Isometrically Equivariant MRE}
The value of an equivariant function $\delta$ on an $\mathcal{X}^n/G$ orbit $[x_0]$ is uniquely determined by the value of $\delta(x_0)$. Informally, the general idea presented below is to find the optimal value of $\delta(x_0) \in \mathcal{X}$ for every possible orbit $[x_0]$ by minimizing the expression $E(L(\theta,\delta)|[x] = [x_0])$.  
 The primary tool used to find  $E(L(\theta,\delta)|[x] = [x_0])$ is the factorization of the base measure $\text{vol}(dx)$ on $\mathcal{X}^n$ into a product of measures on $G$ and $\mathcal{X}^n/G$ \cite{eaton1989group}.

As a first step towards this factorization, the isometry group is given the compact-open topology that is generated by the subbase
\begin{align*}
  V_{K,U} =  \{f \in \text{Iso}(\mathcal{X}): f(K) \subset U\}, \;\; K \subset \mathcal{X}
  \; \text{compact}, \; U \subset \mathcal{X} \; \text{open}.
\end{align*}
 Under this topology the isometry group of any Riemannian manifold is a Lie group  by the Myers-Steenrod theorem \cite{myers1939group}, and thus is a locally compact Hausdorff space. Haar measures therefore exist for $G$. By \cite{bourbaki2004integration} there exists a function $\psi(x) > 0$ such that $\psi(gx) = \Delta(g)\psi(x)$, implying that $\psi(x)\nu(dx) = \text{vol}(dx)$ where $\nu(dx)$ is a relatively invariant measure on $\mathcal{X}^n$ with multiplier $\Delta(g^{-1})$. A factorization of the $\Delta(g)^{-1}$-relatively invariant measure $\nu(dx)$ is given in \cite{eaton1989group} where 
 \begin{align}
     R(\theta,\delta) & = \int_{\mathcal{X}^n} d(\delta(x),EP_\theta)^2 p_n(x|\theta) \text{vol}(dx)  =     \int_{\mathcal{X}^n} d(\delta(x),EP_\theta)^2 p_n(x|\theta)\psi(x)\nu(dx)
     \label{measfact1}
     \\
     & =     \int_{\mathcal{X}^n/G} \bigg(\int_G d(\delta(gx),EP_\theta)^2 p_n(gx|\theta)  \psi(gx)\lambda_R(dg)\bigg) \Tilde{\nu}(d[x]) 
     \label{measfact2}
     \\
     & = \int_{\mathcal{X}^n/G} \Tilde{f}([x]) \Tilde{\nu}(d[x]). 
     \label{measfact3}
 \end{align}
Viewing the inner integral in \eqref{measfact2} as a function $f(x)$ of $x$, the right invariance of $\lambda_R$ shows that $f(hx) = f(x)$ for all $h \in G$. Thus $f(x)$ is constant on $G$-orbits so there exists a function $\Tilde{f}:\mathcal{X}^n/G\rightarrow \mathbb{R}^+$ with $\Tilde{f}([x]) = f(x)$ for all $x \in \mathcal{X}^n$. The factorization Theorem 5.5 of \cite{eaton1989group} asserts the existence of a measure $\Tilde{\nu}$ on $\mathcal{X}^n/G$ such that the integrals in \eqref{measfact1} and \eqref{measfact2} are equal. A regularity assumption is needed for this factorization; the map $(g,x) \rightarrow (gx,x)$ must be proper. The action of the isometry group $\text{Iso}(M)$ of a Riemannian manifold $M$ on $M$ is known to be proper \cite{lee2018introduction}. Moreover, the restriction of a proper map to a closed subset of its domain is also proper. As $G$ is a closed subset of $\text{Iso}(\mathcal{X}^n)$ (see Lemma \ref{grouptoplemma}) it follows that $G$ acts properly on $\mathcal{X}^n$ as is needed for the above factorization. Using the various invariance properties of $d,\delta,EP_\theta$ and $p$ along with the relationships $\Delta(g)\lambda_R(dg) = \lambda_L(dg)$ and $\lambda_L(dg^{-1}) = \lambda_R(dg)$, the expression in \eqref{measfact2} can be rewritten as
\begin{align}
    R(\theta,\delta) & =  \int_{\mathcal{X}^n/G} \bigg(\int_G d(\delta(x),g^{-1}EP_\theta)^2 p_n(x|g^{-1}\theta)  \psi(x) \Delta(g)\lambda_R(dg)\bigg) \Tilde{\nu}(d[x]) 
    \nonumber
    \\
    & = \int_{\mathcal{X}^n/G} \bigg(\int_G d(\delta(x),g^{-1}EP_\theta)^2 p_n(x|g^{-1}\theta)  \psi(x) \lambda_L(dg)\bigg) \Tilde{\nu}(d[x]) 
    \nonumber
    \\
    & = \int_{\mathcal{X}^n/G} \bigg(\psi(x)\int_G d(\delta(x),EP_{g\theta})^2 p_n(x|g\theta)  \lambda_R(dg)\bigg) \Tilde{\nu}(d[x]) .
    \label{measfact4}
\end{align}
As $G$ is transitive over $\Theta$, to find the MRE it suffices to minimizes $R(\theta,\delta)$ at a single value of $\theta$.
The expression in \eqref{measfact4} expresses the risk as a function of the orbit $[x]$, which is enough to determine the explicit form of the MRE. 
\begin{customthm}{1}
\label{MREthm1}
Let $X_1,\ldots,X_n \overset{i.i.d.}{\sim} P_\theta$ be $\mathcal{X}$ valued random objects where $\mathcal{X}$ is a homogeneous Riemannian manifold. Assume that $P_\theta$ lies in the invariant family of distributions $\mathcal{P}$ that is dominated by the Riemannian volume measure on $\mathcal{X}$, and that $\text{Iso}(\mathcal{X})$ acts transitively on $\mathcal{P}$. If the  MRE under the loss function $L(EP_\theta,\delta) = d(EP_\theta,\delta)^2$ exists, it has the form 
 \begin{align}
    \delta_{MRE}(x) = \underset{\delta \in \mathcal{X}}{\argmin} \; E\big(L(\theta,\delta)|[x]\big) = \underset{\delta \in \mathcal{X}}{\argmin} \; \int_G d(\delta,EP_{g\theta})^2 p_n(x|g\theta)  \lambda_R(dg). 
    \label{MREformula}
\end{align}
\end{customthm}

This is a formulation of the classical result that the MRE is the Bayes estimator of $P_\theta$ under a prior distribution for $\Theta$ that is the pushforward of the right Haar measure under the map $g \rightarrow g\theta$ \cite{zidek1969representation,stein1965approximation}. Like the Fr\'echet mean optimization problem, the optimization problem in \eqref{MREformula} is theoretically unwieldy. However, if \eqref{MREformula} has a solution at $x_0$ then it also has a solution for all  $\Tilde{x}$ with $[\Tilde{x}] = [x_0]$ by equivariance.   

The Bayesian setup implied by \eqref{MREformula} has a prior distribution placed on $G$ rather than on $\Theta$. It might be expected that placing a right Haar prior on $G$ is similar to placing a uniform prior over $\Theta$. Due to the transitivity of the action of $G$ on $\Theta$, each left coset of the isotropy group $G_{\theta_0}$ defined by \eqref{isotropydefinition} can be bijectively identified with $\Theta$ by the map $gG_{x_0} \rightarrow gx_0$. Consequently, $\Theta \cong G/G_{\theta_0}$ and $\Theta$ inherits the quotient topology of $G/G_{\theta_0}$ where $\theta_0$ is an arbitrarily chosen point of $\Theta$. A further factorization of the $\lambda_R$ appearing in \eqref{MREformula}, into measures on $G_{\theta_0}$ and $G/G_{\theta_0}$ is possible because $G_{\theta_0}$ is compact. Letting $\lambda_{G_{\theta_0}}$ denote the Haar measure on $G_{\theta_0}$ Corollary 7.4.4 of \cite{wijsman1990invariant} implies
\begin{align}
    \int_G d(\delta & ,EP_{g\theta_0})^2 p(x|g\theta_0)  \lambda_R(dg) \nonumber
    \\
    & =  \int_{G/G_{\theta_0}}\bigg(\int_{G_{\theta_0}} d(\delta,EP_{gh\theta_0})^2p_n(x|gh\theta_0) \lambda_{G_{\theta_0}}(dh)\bigg)\Bar{\nu}(d[g]).
    \label{measfactonG}
\end{align}
As in \eqref{measfact2}, the inner integral in \eqref{measfactonG} viewed as a function of $g$ is constant on the left cosets of $G/G_{\theta_0}$ and so the inner integral is a function of $[g] = gG_{\theta_0}$. The measure $\Bar{\nu}(d[g])$ is the unique $\Delta(g)^{-1}$-invariant measure on $G/G_{\theta_0}$ relative to the action $(a,bG_{\theta_0}) \rightarrow ab G_{\theta_0}$. 

\begin{customcor}{1}
\label{CorMREviaTheta}
Under the same assumptions as Theorem \ref{MREthm1}, if $\Bar{\nu}(d\theta)$ is the unique $\Delta(g)^{-1}$-invariant measure on $\Theta \cong G/G_{\theta_0}$, then the MRE, if it exists, has the form
\begin{align}
    \delta_{MRE}(x) =  \underset{\delta \in \mathcal{X}}{\argmin} \;\; \int_{\Theta} d(\delta,EP_{\theta})^2p_n(x|\theta)\Bar{\nu}(d\theta).
    \label{MREviaTheta}
\end{align}
\end{customcor}
One notable case where \eqref{MREviaTheta} takes a particularly simple form is when $\mathcal{P}$ is parameterized by its Fr\'echet mean, so that $\theta = EP_\theta$ and thus $\Theta = \mathcal{X}$. The induced $G$ action on $\Theta$ is exactly the same as the action of $\text{Iso}(\mathcal{X})$ on $\mathcal{X}$. If in addition $G$ is unimodular so that $\Delta(g)^{-1} = 1$, then the Riemannian volume measure on $\Theta = \mathcal{X}$ is the unique $\Delta(g)^{-1}$-invariant measure under this action and $\hat{\nu}(d\theta) = \text{vol}(d\theta)$.

 Being a Bayes estimator, standard Markov chain Monte Carlo techniques can be utilized to compute the value of the MRE at a given value of $x$. This is done as follows:
 \begin{enumerate}[(1)]
     \item Compute the value of $EP_{\theta_0}$ for some conveniently chosen $\theta_0 \in \Theta$.
     \item Draw a Monte Carlo sample of $g_1,\ldots,g_k$ from the density $p_n(x|g\theta)\lambda_R(dg)$.
     \item Apply each isometry $g_i$ to the point $EP_{\theta_0}$ and obtain the resulting points
     
     $X_i' \coloneqq g_iEP_{\theta_0} = EP_{g_i\theta_0}, \;\; i = 1,\ldots,k$.
     \item Compute the sample Fr\'echet mean of $X_1',\ldots,X_k'$ and take this to be the Monte Carlo approximation of the MRE at $x$, $\delta_{MRE}(x)$.
 \end{enumerate}
 
 The Fr\'echet mean $EP_{\theta_0}$ in step (1) can be found by evaluating the sample Fr\'echet mean of a large number of i.i.d.\@ Monte Carlo draws from $P_{\theta_0}$. 
 If $\Theta \cong \mathcal{X}$ and $G$ is unimodular then steps (2) and (3) in the above procedure can be replaced by:
  \begin{enumerate}
     \item[($2'$)] Draw a Monte Carlo sample of $\theta_1,\ldots,\theta_k$ from the density $p_n(x|\theta)\text{vol}(d\theta)$.
     \item[($3'$)] Compute the Fr\'echet means $EP_{\theta_i}, \; i = 1,\ldots,k$ and take $X_i' \coloneqq EP_{\theta_i}$.
 \end{enumerate}
However, it is often easier to work in $G$ since in $G$ only $EP_{\theta_0}$ needs to be calculated to find the MRE, while in $\Theta$ each of $EP_{\theta_1},\ldots,EP_{\theta_k}$ must be computed. The particular choice of Monte Carlo algorithm used in step $(2)$ or $(2')$ is problem dependent. When $n$ is small and $G$ is compact a simple method for obtaining $g_1,\ldots,g_k$ in $(2)$ is to use uniform $\lambda_R$ proposals in a Metropolis-Hastings algorithm. As $n$ grows larger the posterior $p_n(x|g\theta)$ becomes more peaked and more sophisticated proposals are needed.

The invariant estimation problem formulated here and its solution easily generalize to the case of estimating a $k$-Fr\'echet mean for an arbitrary $k$. The loss function used can also any positive power of $d(\cdot,\cdot)$. Moreover, the MRE depends on the choice of the Riemannian distance $d$ only through the isometry group $\text{Iso}_d(\mathcal{X})$ and the objective function of the optimization problem \eqref{MREformula}. If $d$ is a Riemannian distance on $\mathcal{X}$ and $\Tilde{d}$ is another, not necessarily Riemannian, distance on $\mathcal{X}$ with $\text{Iso}_{\Tilde{d}}(\mathcal{X})$ a closed subgroup of $\text{Iso}_{d}(\mathcal{X})$, then $\text{Iso}_{\Tilde{d}}(\mathcal{X})$ acts properly on $\mathcal{X}$ and the above factorizations of $\text{vol}_{d}(dx)$ remain valid. In such cases, the MRE with a loss function $\Tilde{d}(\delta,EP_\theta)^2$ is given by \eqref{MREformula} with $d$ replaced by $\Tilde{d}$ and with $\lambda_R$ replaced by the right Haar measure $\Tilde{\lambda}_R$ on $\text{Iso}_{\Tilde{d}}(\mathcal{X})$. It is up to the statistician to choose an appropriate distance function that reflects the loss for the problem at hand.

\subsection{Applications}
\subsubsection{von Mises-Fisher Distributions on the Sphere and Hyperbolic Space}
One of the simplest examples of a non-Euclidean $k$-dimensional homogeneous Riemannian manifold is the $k$-sphere, $\mathbb{S}^k \subset \mathbb{R}^{k+1}$. The distance between two-points on $\mathbb{S}^k$ is defined to be the length of the shortest path between these points that lies on the sphere. Any orthogonal transformation in $\text{O}(k+1)$ maps $\mathbb{S}^k$ to itself and preserves the lengths of paths on the sphere.  In fact $\text{Iso}(\mathbb{S}^k) = \text{O}(k+1)$, which agrees with intuition for $\mathbb{S}^2$ since reflections or rotations of a sphere do not distort the geometry of the sphere. It is also clear that $\text{Iso}(\mathbb{S}^k)$ is transitive since a point can always be rotated to any other point on the sphere. The von Mises-Fisher family of distributions on the sphere have densities $p(x|\mu,\kappa) \propto \exp(\kappa x^\intercal \mu)$ parameterized by  $(\mu,\kappa) \in \Theta \coloneqq \mathbb{S}^k \times \mathbb{R}^+ $  with respect to the volume measure. The parameters $\mu$ and $\kappa$ are interpreted as location and concentration parameters respectively. If $U \in \text{O}(k+1)$ is an isometry of $\mathbb{S}^k$ then $UP_{(\mu,\kappa)}$ has the density $p(U^{-1}x|\mu,\kappa) \propto \exp(\kappa (U^\intercal x)^\intercal \mu)$ so that $UP_{(\mu,\kappa)} = P_{(U\mu,\kappa)}$. The set of orbits for the entire von Mises-Fisher family is $\Theta/G \cong \{\kappa: \kappa > 0\} = (0,\infty)$ where $G$ is transitive over any subfamily of distributions that have a fixed value of $\kappa$.

As might be expected, the Fr\'echet mean of $P_{(\mu,\kappa)}$ is $\mu$ for $\kappa > 0$. The sphere is a two-point homogeneous space, meaning that for any $p_1,q_1,p_2,q_2 \in \mathbb{S}^n$ with $d(p_1,q_1) = d(p_2,q_2)$ there exists an isometry taking $p_1$ to $p_2$ and $q_1$ to $q_2$. The following theorem provides a way to determine the Fr\'echet mean of a specific class of distributions on two-point homogeneous spaces. This theorem extends results for shape spaces that appear in \cite{kendall2009shape}.
\begin{customthm}{2}
If $\mathcal{X}$ is a two-point homogeneous Riemannian manifold and $X$ is a random object taking values in $\mathcal{X}$ with the density $f\big(d(x,\mu)\big)h(x)$ with respect to the Riemannian volume measure, then $EX = \mu$ if $f$ is a decreasing function. 
\end{customthm}
\begin{proof}
See the Appendix for a proof of Theorem \ref{generalFrecmnthm}. 
\end{proof}
As the distance between points $x,y \in \mathbb{S}^k$ is given by $d(x,\mu) = \text{arccos}(\langle x, \mu \rangle)$, the von Mises-Fisher density can be expressed as $p(x|\mu,\kappa) \propto \exp(\kappa \cos(d(x,\mu))$ and Theorem \ref{generalFrecmnthm} implies that $EP_{(\mu,\kappa)} = \mu$.

For a fixed value of $\kappa = \kappa_0$ and an i.i.d.\@ sample $X_1,\ldots,X_n \overset{i.i.d.}{\sim} P_{(\mu,\kappa_0)}$, the MRE of $EP_{(\mu,\kappa_0)} = \mu$ under the squared distance loss can be found using \eqref{MREviaTheta}: 
\begin{align}
    \delta_{MRE}(X) =  \underset{\delta \in \mathbb{S}^n}{\text{argmin}} \;\; \int_{\mathbb{S}^n} d(\delta,\mu)^2 \exp\bigg(\kappa_0\Vert S_n \Vert \big(\mu^\intercal \frac{S_n}{\Vert S_n \Vert}\big)\bigg) \text{vol}(d\mu),
    \label{vonmisMRE}
\end{align}
 where $S_n = \sum_{i = 1}^n X_i$. Note that the compactness of $\text{O}(k+1)$ implies it is unimodular and the remarks that immediately follow Corollary \ref{CorMREviaTheta} apply here. Due to the conjugacy of the von Mises-Fisher distribution, $\delta_{MRE}(X)$ can be recognized as the Fr\'echet mean of $P_{(S_n/\Vert S_n \Vert,\kappa_0\Vert S_n \Vert)}$ conditional on $S_n$, which is $S_n/\Vert S_n\Vert$. The MRE in this case is equal to the MLE of $\mu$ and moreover does not depend on the orbit $\kappa_0$. Now suppose that instead of $d$, the extrinsic distance $\Tilde{d}(x,y) = \Vert x - y \Vert$ is used, where $\mathbb{S}^k$ is viewed as an embedded submanifold of $\mathbb{R}^{k+1}$ \cite{bhattacharya2003large}. As $\text{Iso}_{\Tilde{d}}(\mathbb{S}^k) = \text{O}(k+1) = \text{Iso}_d(\mathbb{S}^k)$ the MRE is given by $\eqref{vonmisMRE}$ with $d$ replaced by $\Tilde{d}$. The MRE is identical to the previous case as the Fr\'echet mean of $P_{(S_n/ \vert S_n \Vert,\kappa_0 \Vert S_n \Vert)}$ with respect to $\Tilde{d}$ is also $S_n/\Vert S_n \Vert$. A similar result can be found in \cite{sengupta1998best} where the simultaneous estimation problem of $\mu_i \in \mathbb{S}^1, \; i = 1,\ldots,m$ is considered given independent observations of $X_{ij} \sim P_{(\mu_i,\kappa)}, i = 1,\ldots,m, \; j = 1,\ldots,n_i$. It is shown that the MRE of $(\mu_1,\ldots,\mu_m)$ is the maximum likelihood estimator in this setting. Our result generalizes the solution of \cite{sengupta1998best} in the $m = 1$ case to higher dimensional spheres and has a coordinate free derivation. 

The condition of two-point homogeneity implies that $\mathcal{X}$ is a symmetric space \cite{varma1965two}. Other examples of symmetric, two-point homogeneous spaces include Lie groups with bi-invariant Riemannian metrics, Grassmannians and hyperbolic space \cite{lee2018introduction}. Hyperbolic space is of special interest in differential geometry as it is the  ``model space" of negative curvature. The construction of the hyperboloid model of hyperbolic space proceeds identically to the construction of spheres except that the  Minkowski pseudo-inner product, $(x,y) = \sum_{i = 1}^kx_iy_i - x_{k+1}y_{k+1}$, is used instead of the Euclidean inner product. The hyperboloid is defined as the collection of points $\mathbb{H}^k(R) = \{x \in \mathbb{R}^{k+1}: (x,x) = -R^2, \; x_{k+1} > 0\}$.  Distances in the hyperboloid model are given by $d(x,y) = R \: \text{arcosh}(-(x,y)/R^2)$, a formula reminiscent of the angular distance between points on a sphere  $d_{\mathbb{S}^k(R)}(x,y) = R\,\arccos(\langle x, y \rangle/R^2)$. The von Mises-Fisher analogue on the hyperboloid is the hyperbolic distribution with density $p(x|\kappa,\mu) \propto \exp(\kappa \: (x,\mu)/R^2), \; \mu \in \mathbb{H}^k(R), \; \kappa > 0$ with respect to $\text{vol}(dx)$ \cite{jensen1981hyperboloid, barndorff1978hyperbolic}.  An application of Theorem \ref{generalFrecmnthm} shows that the Fr\'echet mean of a hyperbolic distribution is $EP_{(\kappa,\mu)} = \mu$.  The connected component containing the identity of the isometry group $\text{Iso}(\mathbb{H}^k) = \text{O}^+(k,1)$ is transitive over $\mathbb{H}^k$ and thus is transitive over $\Theta$ for a fixed value of $\kappa$.  It then suffices to consider only the action of this connected component subgroup $\text{SO}^+(k,1)$ when finding the MRE. Despite it not being compact, $\text{SO}^+(k,1)$ is unimodular because it is a semisimple Lie group \cite{knapp2013lie}. By \eqref{MREviaTheta} and the same argument as in the spherical case, if $S_n = \sum_{i = 1}^n X_i$ the MRE for $\mu$ is  $RS_n/\big(-(S_n,S_n)\big)^{1/2}$, which also coincides with the MLE.

\subsubsection{Langevin Distribution on the Stiefel Manifold}
The equivariance of the Fr\'echet mean functional also provides another tool that can be used to determine the Fr\'echet mean of a distribution.
\begin{customlemma}{2}
\label{IsotropFrecmnlemma}
The isotropy group $G_\theta$ is contained in the isotropy group $G_{EP_\theta}$. If $EP_\theta$ is not a singleton set this containment still applies where $G_{EP_\theta}$ is the isotropy group with respect to the action of $G$ on the power set of $\mathcal{X}$.
\end{customlemma}  
\begin{proof}
If $g \in G_{\theta}$ then $   gEP_{\theta} = EgP_{\theta} = EP_{g\theta} = EP_{\theta}$ so $g \in G_{EP_{\theta}}$. 
\end{proof}

The above lemma is most useful when $G_\theta$ is large, thereby limiting the possible values of $EP_\theta$. As an application, consider the Stiefel manifold $\text{V}_k(\mathbb{R}^p)$ that consists of all matrices $X \in \mathbb{R}^{p \times k}$ with $X^\intercal X = \text{I}_k$. Viewing the Stiefel manifold as an embedded submanifold of $\mathbb{R}^{p \times k}$ it can inherit either the extrinsic Frobenius norm distance or the intrinsic, induced Riemannian distance, although other distances are also useful \cite{edelman1998geometry}. As the Frobenius norm satisfies $\Vert UAV^\intercal - UBV^\intercal \Vert_F = \Vert A - B \Vert_F$ for all $A,B \in \mathbb{R}^{p \times k}, \; U \in \text{O}(p), \; V \in \text{O}(k)$, the Stiefel manifold contains the direct product of $\text{O}(p)$ and $\text{O}(k)$ in its isometry group. Similarly, the isometry group under the induced Riemannian distance also contains $G \coloneqq \text{O}(p) \times \text{O}(k)$.

The Langevin distribution on $\text{V}_k(\mathbb{R}^p)$ generalizes the von Mises-Fisher distribution and has a density $p(x|\theta) \propto \text{etr}(x^\intercal \theta), \; \theta \in \mathbb{R}^{p \times k}$ with respect to $\text{vol}(dx)$. Under the action of $(U,V) \in \text{O}(p) \times \text{O}(k)$, the Langevin distribution transforms according to $(U,V)P_\theta = P_{U\theta V^\intercal}$. A Langevin subfamily of interest is $\theta = \lambda H, \; H \in \text{V}_k(\mathbb{R}^p), \; \lambda > 0$. In this family each column vector of $X$ is concentrated around the corresponding column vector of $H$ and all columns of $X$ concentrate around the columns of $H$ by equal amounts. Notice that if $H_0 = [\text{I}_k,0_{p-k,k}]^\intercal$ then $G_{\lambda H_0}$ contains matrices of one of the two forms
\begin{align*}
   \bigg( \begin{bmatrix}
    \text{I}_k & 0 
    \\
    0 & V
    \end{bmatrix}, \text{I}_{k}\bigg), \; V \in \text{O}(p - k),\;\; \text{or} \;\;\; (T_{ij},T_{ij}), \; 1 \leq i,j \leq k,
\end{align*}
where the $T_{ij} = \text{I} - e_ie_i^\intercal - e_je_j^\intercal + e_ie_j^\intercal + e_je_i^\intercal$ are transposition matrices that permute the rows and columns of a matrix under left and right multiplication by $T_{ij}$. It is seen that the only elements of $\text{V}_k(\mathbb{R}^p)$ that are fixed by these group elements are matrices of the form $[\text{diag}(\pm 1), 0_{p-k,k}]^\intercal$. If unique, $EP_{\lambda H_0}$ must be one of these matrices by Lemma \ref{IsotropFrecmnlemma}. By considering left multiplication by reflection matrices $\text{diag}(1,\ldots,-1,\ldots,1)$, it is clear that out of these matrices, $H_0 = [\text{diag}(1), 0_{p-k,k}]^\intercal$ is the closest on average to $X$ so $EP_{\lambda H_0} = H_0$ for all $\lambda > 0$. By the equivariance of $E(\cdot)$, $EP_{\lambda H} = H$ for any $H \in \text{V}_k(\mathbb{R}^p)$.

If either the squared intrinsic or extrinsic distances are used as loss functions, then given a single observation $X \sim P_{\lambda H}$, the MRE of $H$ is $X$. This follows by the conjugacy of this family of Langevin distributions, whose density functions are symmetric in $X$ and $H$. By the remarks following \eqref{MREviaTheta}, the MRE of $EP_{\lambda H}$ is the Fr\'echet mean $EP_{\lambda X} = X$, where again we note that $G$ is unimodular. An implication of this result is that $X$ is an admissible estimator of $H$ under these losses, because it is a Bayes estimator under a proper prior distribution. That the prior is proper is a direct consequence of the finiteness of the right Haar measure on the compact group $G$. When $n > 1$ this particular subfamily of distributions is no longer conjugate and the posterior Fr\'echet mean implied by \eqref{MREviaTheta} must instead be found numerically.  

\section{Estimation Under a Non-transitive Action}
\label{SecNontrans}

\subsection{Adaptive Equivariant Estimator}
\label{AdaptiveeqSubsec}
Many models $\mathcal{P}$ of interest have a $G$-action that does not act transitively on $\mathcal{P}$. Recall that the risk function of any equivariant estimator $\delta$ is constant over $G$-orbits of the parameter space, $R(g\theta,\delta) = R(\theta,\delta)$ for all $g \in G$. Transitivity of $G$ ensures that the risk functions of equivariant estimators can be totally ordered. In the non-transitive setting the risk functions of equivariant estimators have a total ordering when restricted to a single orbit. That is, if $G$ acts on the family $\mathcal{P}$ then for any $\theta_0 \in \Theta$, $G$ will act transitively on the subfamily $\mathcal{P}_{[\theta_0]} \coloneqq \{P_{\theta}:[\theta] = [\theta_0]\} \subset \mathcal{P}$ and the results from the previous section directly apply to $\mathcal{P}_{[\theta_0]}$.

This suggests the possibility of using a two-step estimation procedure where the $\Theta$-orbit that contains the true $\theta$ is first estimated and then the MRE is computed conditional on this estimated orbit. Formalizing this, for each $\theta_0 \in \Theta$ let $\delta_{[\theta_0]}$ be the MRE for the sub-family of distributions $\mathcal{P}_{[\theta_0]}$ where it is assumed that this MRE exists for every such sub-family $\mathcal{P}_{[\theta_0]}$. Let $\widehat{[\theta]}:\mathcal{X}^n \rightarrow \Theta/G$ be an estimator of the true orbit of $\theta$. We define the $\widehat{[\theta]}$-adaptive MRE to be the estimator $\delta_{\widehat{[\theta]}}$. The orbit estimators that we use in the following sections are derived from estimators $\hat{\theta}(X_1,\ldots,X_n)$ of the full parameter $\theta$ where $\widehat{[\theta]}(X_1,\ldots,X_n) = [\hat{\theta}(X_1,\ldots,X_n)]$.

\begin{customlemma}{3}
If $\widehat{[\theta]}$ is $G$-invariant, the adaptive MRE $\delta_{\widehat{[\theta]}}$ is equivariant. In particular, if $\widehat{[\theta]}(X) = [\hat{\theta}(X)]$ for an equivariant estimator $\hat{\theta}:\mathcal{X}^n \rightarrow \Theta$ then the adaptive MRE is equivariant. 
\end{customlemma}
\begin{proof}
If $\hat{\theta}:\mathcal{X}^n \rightarrow \Theta$ is equivariant then $\widehat{[\theta]} = [\hat{\theta}]$ is invariant since
\begin{align*}
    \widehat{[\theta]}(gX) = [\hat{\theta}(gX)] = [g\hat{\theta}(X)] = [\hat{\theta}(X)] = \widehat{[\theta]}(X).
\end{align*}
  For an invariant $\widehat{[\theta]}$ the adaptive MRE is equivariant as
\begin{align*}
\delta_{\widehat{[\theta]}(gX)}(gX) = \delta_{\widehat{[\theta]}(X)}(gX) = g\delta_{\widehat{[\theta]}(X)}(X).
\end{align*}
\end{proof}
It is expected that if the orbit estimate is accurate, the adaptive MRE will perform similarly to the MRE under the sub-family $\mathcal{P}_{[\theta]}$ containing the true value of $\theta$. Thus, it is hypothesized that $\delta_{\widehat{[\theta]}}$ will perform well when the sample size is large.  

The estimation procedure described above is closely related to empirical Bayes estimation. The orbit $[\theta]$ can be viewed as a hyperparameter for the Bayesian model $X_1,\ldots,X_n \overset{i.i.d.}{\sim} P_{g\theta}, \; g \sim \lambda_R$. As the Bayes estimator of $EP_\theta$ for a fixed orbit $[\theta]$ is $\delta_{[\theta]}$, the empirical MRE $\delta_{\widehat{[\theta]}}$ can be viewed as a empirical Bayes estimate of $EP_\theta$ that uses the data to estimate the hyperparameter $[\theta]$. Connections between equivariant estimation under a non-transitive action and the James-Stein estimator are discussed in \cite{beran1996stein}. There the action of $\text{O}(p)$ on $\mathbb{R}^p$ is considered when estimating $\mu$ given a single observation from the model $X \sim N_p(\mu,\text{I})$. If $X \sim N_p(\mu, \text{I})$ then $UX \sim N_p(U\mu,\text{I})$ for $ U \in \text{O}(p)$ and the $\mathcal{P}$ orbits of this model can be indexed by $[\mu] = \Vert \mu \Vert$. For a fixed orbit $[\mu ] = \Vert \mu \Vert$, let $\delta_{\Vert \mu \Vert}$ be the MRE  restricted to take values in the sphere of radius $\Vert \mu \Vert$. That is, $\delta_{\Vert \mu\Vert}$ is MRE out of all estimators whose action space is $\mathbb{S}^p(\Vert \mu \Vert)$. It is shown in \cite{beran1996stein} that the James-Stein estimator is equal to $\delta_{\widehat{\Vert \mu \Vert}}$ when $\widehat{\Vert \mu \Vert} = \Vert X \Vert^2 - p  + 2$. Note that $X$ is the MRE under the full, transitive isometry group of $\mathbb{R}^n$ which includes the translations $x \rightarrow x + a$. As the James-Stein estimator dominates $X$ for $p > 2$, this demonstrates how it can be beneficial to minimize the size of the group under which equivariance of an estimator is required.

\subsection{Numerical Illustration: Positive Definite Matrices}

Covariance estimation of $p \times p$ positive-definite matrices $\mathcal{S}^p_+$ is another classical setting where non-transitive group actions are useful. Suppose $X \sim \text{Wishart}_p(n,\Sigma)$ with $n$ known and it is desired to estimate $\Sigma$ under an orthogonally invariant loss $L(U \Sigma U^\intercal,U \delta U^\intercal) = L(\Sigma, \delta), \; U \in \text{O}(p)$. The orthogonal group acts on a covariance matrix $X$ by $X \rightarrow UXU^\intercal$ where $UXU^\intercal \sim \text{Wishart}_p(n,U\Sigma U^\intercal)$. Each parameter space orbit can be indexed by the eigenvalues of $\Sigma$ so that $\Theta/G = \mathbb{R}^p_+$. It is possible to extend this action from $\text{O}(p)$ to $\text{GL}_p(\mathbb{R})$ which makes the resulting action transitive over $\Theta$. However, not all loss functions will be invariant under the full $\text{GL}_p(\mathbb{R})$ action.

In this subsection we consider the loss function given by the log-Euclidean distance on $\mathcal{S}_+^p$ defined by $d(X,Y) \coloneqq \Vert \log(X) - \log(Y) \Vert_F$ where $\log(\cdot)$ is the matrix logarithm and $\Vert \cdot \Vert_F$ is the Frobenius norm. The Euclidean distance applied to covariance matrices can exhibit a swelling effect where it is possible that the determinant of a mean of matrices is larger than the determinant of any individual matrix in the mean. Motivation for using the log-Euclidean distance stems in part as a way to mitigate this swelling effect. The log-Euclidean distance is especially useful in medical imaging applications where determinants of covariance matrices have direct physical interpretations.

The matrix logarithm is a bijection between $\mathcal{S}_+^p$ and the set of $p \times p$ symmetric matrices. If $X$ has the eigendecomposition $X = U\text{diag}(\lambda_1,\ldots,\lambda_p)U^\intercal$ the matrix logarithm is defined by $\log(X) = U \text{diag}(\log(\lambda_1),\ldots,\log(\lambda_p))U^\intercal$. It follows that $\log(\cdot)$ is $\text{O}(p)$ equivariant, $\log(UXU^\intercal) = U\log(X)U^\intercal$, which shows that $\text{O}(p)$ is contained in $\text{Iso}(\mathcal{S}_+^p)$ for the log-Euclidean metric. The log-Euclidean metric is not fully $\text{GL}_p(\mathbb{R})$ invariant, contrasting with other commonly used losses such as  $L_1(\Sigma,\delta) = \text{tr}(\delta\Sigma^{-1}) - \log(\vert \delta \Sigma^{-1}\vert)$ and $L_2(\Sigma,\delta) = \text{tr}(\delta\Sigma^{-1} -  \text{I})$. Working in log coordinates, the log-Euclidean Fr\'echet mean $EX$ of $X$ corresponds to the Euclidean mean of $\log(X)$, $\log(EX) = E_{Euc}(\log(X))$.

  The decision problem considered in the simulation study below is to estimate the Fr\'echet mean of $X/n$ under the log-Euclidean loss $L(\Sigma,\delta) = d(\Sigma,\delta)^2$  when $X \sim \text{Wishart}_p(n,\Sigma)$. As $d$ is induced from a Riemannian metric on $\mathcal{S}^p_+$ \cite{arsigny2006log} the theory from Sections \ref{SecEquivariantEst} and \ref{AdaptiveeqSubsec} applies and an adaptive equivariant estimator can be employed. Alternative estimators of $E(X/n)$ include the sample Fr\'echet mean of $X/n$, which is simply $X/n$, and the MLE, which is $E(Y/n)$ where $Y \sim \text{Wishart}_p(n,X/n)$ conditionally on $X$. The eigenvalues of any equivariant estimate of $\Sigma$, such as the MLE $X/n$, can serve as an orbit estimate for the adaptive MRE. An orbit estimate can also be obtained from the sample Fr\'echet mean by a method of moments procedure where the method of moments estimator for $\Sigma$ satisfies the equation
\begin{align*}
    X/n = E(Y/n), \;\; Y \sim \text{Wishart}_p(n,\Sigma_{MoM}).
\end{align*}
Simulation results show that this equation has an approximate solution for $\Sigma_{MoM}$. The resulting adaptive MRE is given by $\delta_{[\Sigma_{MoM}]}$, where we recall from Subsection \ref{AdaptiveeqSubsec} that $\delta_{[\Sigma_0]}$ is the MRE of the sub-model where $\Sigma$ is restricted to lie in the same $\text{O}(p)$ orbit as $\Sigma_0$.

The estimated risks of the sample Fr\'echet mean, MLE and adaptive MRE with MLE and $\Sigma_{MoM}$ orbit estimates are shown in Table \ref{SPDRiskTable} for $\Sigma = \text{diag}(1,\ldots,p)$. These risks were computed by averaging the observed losses of each estimator over $500$ different simulated data sets of sample size $n \in \{5,10,40\}$. Given an orbit, the MREs are computed by the procedure described in Section \ref{SecEquivariantEst} where a Metropolis Hastings algorithm over $G = \text{O}(p)$ with uniform proposals is run for $1500$ iterations. In each case it is seen that the adaptive MRE improves upon the estimator from which its orbit is derived from. The adaptive MRE with the $\Sigma_{MoM}$ orbit performs especially well under every scenario. Surprisingly, even in the $n = 5$ case where the orbit estimates may be inaccurate the performance of the adaptive MREs is superior to the sample Fr\'echet mean and MLE.

\begin{table}
    \centering
    \begin{tabular}{|c|c|c|c|c|}
    \hline
  &  \multicolumn{4}{|c|}{$\delta$}
    \\
    \hline
    & Sample Fr\'echet & MLE & MRE, $\widehat{[\theta]} = X/n$&  MRE, $\widehat{[\theta]} = \Sigma_{MoM}$
    \\
    \hline
    \hline
    $p = 2$, $n = 5$ & 1.796 & 2.234 & 2.004 & 0.216
\\
\hline
  $p = 2$, $n = 10$ & 0.723 & 0.803 & 0.715 & 0.168
\\
\hline
  $p = 2$, $n = 40$ & 0.143 & 0.146 & 0.144 & 0.055    
\\
\hline
  $p = 4$, $n = 5$ & 9.754 & 16.021 & 14.555 & 1.379 
\\
\hline
  $p = 4$, $n = 10$ & 2.664 &  3.361 & 2.795 & 0.890 
\\
\hline
  $p = 4$, $n = 40$ & 0.513 & 0.548 & 0.506 & 0.289 
\\
\hline
    \end{tabular}
    \caption{Risk of estimator $\delta$, $R(\theta,\delta)$}
    \label{SPDRiskTable}
\end{table}

\subsection{Numerical Illustration: Data on the Torus}

A product Riemannian manifold $\mathcal{X} = M_1 \times \cdots \times M_p$ is the Cartesian product of the Riemannian manifolds $M_i$ and is equipped with the distance function \begin{align*}
    d_{\mathcal{X}}((x_1,\ldots,x_p),(y_1,\ldots,y_p))^2 \coloneqq \sum_{i = 1}^p d_{M_i}^2(x_i,y_i).
\end{align*}
Examples of commonly used product manifolds include $\mathbb{R}^p = \mathbb{R} \times \cdots \times  \mathbb{R}$ and the (flat) $p$-torus, $\mathbb{T}^p \coloneqq \mathbb{S}^1 \times \cdots \times \mathbb{S}^1$ which is the primary focus of this section. Tori have been used to represent multivariate angle measurements such as torsion angles in proteins and other biological molecules \cite{boomsma2008generative}.

In a general product manifold $\mathcal{X} = M_1 \times \cdots \times M_p$ the Fr\'echet mean of a random object $X = (X^{(1)},\ldots,X^{(p)}) \in \mathcal{X}$ behaves just like the mean of a random vector. The mean of a product is the product of the marginal means, $E_{d_{\mathcal{X}}}X = (E_{d_{M_1}}X^{(1)},\ldots,E_{d_{M_k}}X^{(p)})$. In particular, sample Fr\'echet means can be computed by computing the marginal sample Fr\'echet means. If each $M_i$ is a homogeneous space then $\mathcal{X}$ is also homogeneous since $\text{Iso}(M_1) \times \cdots \times \text{Iso}(M_p) \subset \text{Iso}(\mathcal{X})$ acts transitively on $\mathcal{X}$. To ease notation we denote $\text{Iso}(M_1) \times \cdots \times \text{Iso}(M_p)$ by $\Tilde{G}$. Equivariant estimation of the Fr\'echet mean $EX$ in the product manifold setting is most interesting when the components of $X$ are dependent. If the $X^{(i)}$ are independent and each $\text{Iso}(M_i)$ acts transitively over the family of marginal distributions of $X^{(i)}$, equivariant estimation of $EX$ on $\mathcal{X}$ is no better than separately performing equivariant estimation on the marginal distributions.

\begin{figure}[h]
    \centering
    \includegraphics[scale = 0.95]{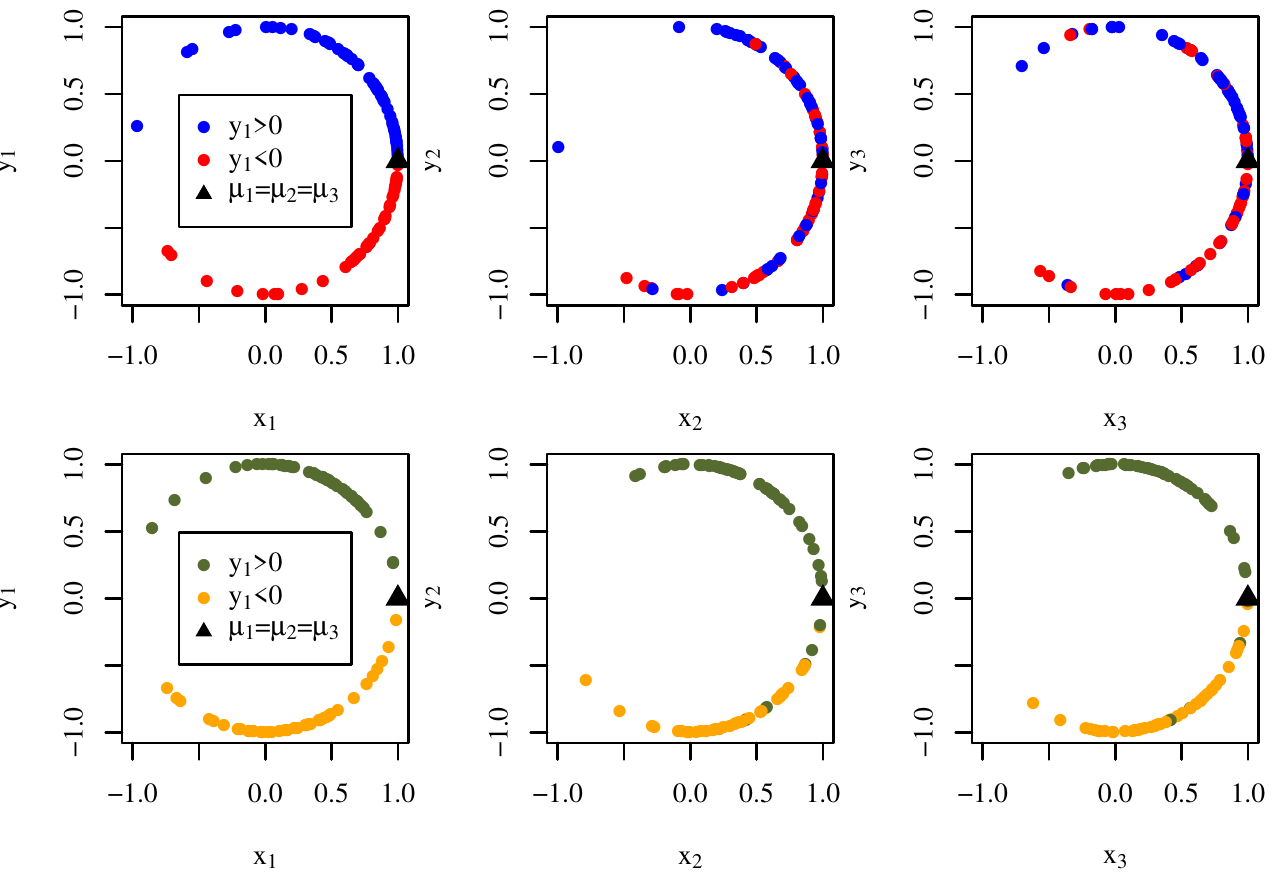}
    \caption{First row: $\kappa_i = 2, \; \lambda_{ij} = 1$, Second row: $\kappa_i = 2, \; \lambda_{ij} = 3$}
    \label{TorusSimsFig}
\end{figure}

The decision problem of interest in this section is to estimate the Fr\'echet mean of a random object $X$ taking values in the $p$-torus under the squared distance loss $d_{\mathbb{T}^p}(EX,\delta)^2$. The distribution of $X$ is assumed to have a density with respect to $\text{vol}(dx)$ on $\mathbb{T}^p$ of
\begin{align}
\label{torusdistribution}
    &p(x|\mu, \kappa,\Lambda)  \propto  \exp\bigg( \sum_{i = 1}^p \kappa_i x_i^\intercal \mu_i +  \sum_{i < j}\lambda_{ij} x_i^\intercal (R\mu_i(R\mu_j)^\intercal)x_j\bigg),\;\;\;  R  =  \begin{bmatrix} 0 & -1
    \\
    1 & 0\end{bmatrix},
    \\
    & x_i,\mu_i \in \mathbb{S}^1, \; \kappa_i > 0,\; \lambda_{ij} \in \mathbb{R}, \; \Theta = \{(\mu,\kappa,\Lambda):\mu \in \mathbb{T}^p, \kappa \in \mathbb{R}_+^{p},\Lambda \in \mathbb{R}^{p(p-1)/2} \}.
    \nonumber
\end{align}
This family of distributions is a multivariate generalization of the von Mises-Fisher distribution on the $p$-torus that also induces dependence between the $\mathbb{S}^1$ components of $\mathbb{T}^p$ \cite{singh2002probabilistic, mardia2008multivariate}. The first term in the exponent in \eqref{torusdistribution} is equal to $\sum_i \kappa_i \cos(\angle x_i - \angle \mu_i)$ while the second term is $\sum_{i < j}\lambda_{ij}\sin(\angle x_i - \angle \mu_i)\sin(\angle x_j  - \angle \mu_j)$, where $\angle$ gives the angle of a point in $\mathbb{S}^1$ relative to $(1,0)$.  The first term reflects the affinity of $x_i$ to be close to $\mu_i$ and the second term controls the correlation between $\angle x_i - \angle \mu_i$ and $\angle x_j - \angle \mu_j$. Figure \ref{TorusSimsFig} illustrates how the value of the $\lambda_{ij}$'s alters this distribution for $p = 3$. It is seen that when the $\lambda_{ij}$'s have the same sign and are large relative to the $\kappa_i$'s. observations in each respective circle tend to cluster more around $\pm R\mu_i$ while if the $\lambda_{ij}$'s are smaller then observations cluster around the $\mu_i$'s. Fr\'echet means for this distribution are therefore most meaningful when the magnitude of the $\lambda_{ij}$'s is not too large.

An isometry $g = g_1 \times \cdots \times  g_p \in \Tilde{G}$ acts on this family by $gP_{(\mu,\kappa,\Lambda)} = P_{(g\mu, \kappa,\Lambda)}$, implying that $\Theta/\Tilde{G} = \mathbb{R}^p_+ \times \mathbb{R}^{p(p-1)/2}$. The normalizing constant of \eqref{torusdistribution} is not known, however it does not depend on $\mu$ since the normalizing constants associated with $P_{(\mu,\kappa,\Lambda)}$ and $gP_{(\mu,\kappa,\Lambda)}$ must be the same. The maximum likelihood estimate of $\theta = (\mu,\kappa,\Lambda)$ can be found by first finding the restricted maximum likelihood estimate $\hat{\mu}_{\kappa,\Lambda} \coloneqq \text{argmax}_{\mu} \ell(\mu,\kappa,\Lambda)$ using the known $\mu$-gradients of the likelihood. The full maximum likeihood estimates can then be found numerically by maximizing the profile likelihood $\ell(\hat{\mu}_{\kappa,\Lambda},\kappa,\Lambda)$ over $\kappa$ and $\Lambda$. 

\begin{table}
    \centering
    \begin{tabular}{|c|c|c|}
    \hline
    Parameter Values $\theta$ &  \multicolumn{2}{|c|}{$\delta$}
    \\
    \hline
    & Sample Fr\'echet & MLE
    \\
    \hline
        \hline
    $\kappa_i = 2$, $\lambda_{ij} = 1$, $n = 25$ & 1.73 & 1.44
\\
$\kappa_i = 2$, $\lambda_{ij} = 1$, $n = 5$ & 1.88 & 1.82
\\
$\kappa_i = 2$, $\lambda_{ij} = 3$, $n = 25$ & 8.98 & 1.33
\\
$\kappa_i = 2$, $\lambda_{ij} = 3$, $n = 15$ & 8.14 & 0.99
\\
$\kappa_i = 2$, $\lambda_{ij} = 3$, $n = 5$ & 3.04 & 0.86
\\
\hline
    \end{tabular}
    \caption{Ratio of risks of estimator $\delta$ to the adaptive MRE, $R(\theta,\delta)/R(\theta,\delta_{[(\hat{\kappa},\hat{\lambda})]})$}
    \label{TorusRiskTable}
\end{table}

By estimating the orbit by the MLE $(\hat{\kappa},\hat{\lambda})$, the adaptive MRE procedure is applicable here. Notice that the conditional distribution $p(x_i|x_{-i},\mu,\kappa,\lambda)$ is a von Mises-Fisher distribution so that Monte Carlo samples from $P_{(\mu,\kappa,\lambda)}$ can be efficiently drawn using Gibbs sampling \cite{mardia2008multivariate}. Moreover, the posterior $p(\mu|x,\kappa,\lambda)$ under a uniform prior density is symmetric in $x$ and $\mu$, implying that Gibbs sampling can also be used to sample from the posterior distribution in \eqref{MREviaTheta}. Given Monte Carlo samples of $\mu_1,\ldots,\mu_k$ from this posterior distribution the adaptive MRE is the sample Fr\'echet mean of $EP_{(\mu_1,\hat{\kappa},\hat{\lambda})},\ldots,EP_{(\mu_k,\hat{\kappa},\hat{\lambda})}$ where $\hat{\kappa},\hat{\lambda}$ are the MLE orbit estimates of $\kappa$ and $\lambda$. Table \ref{TorusRiskTable} provides simulation results for $n \in  \{5,25\}$, $p = 3$ and $\kappa_i$ and $\lambda_{ij}$-values that are constant across $i$ and $j$.  The risks for each estimator are estimated by averaging the observed loss across $1000$ simulated data sets from \eqref{torusdistribution}. Without loss of generality $\mu_i$ is taken to be $(0,1)$ for all $i$ in the simulations, since the risks of the estimators are independent of $\mu$ by equivariance. It is seen that in most cases the adaptive MRE performs significantly better than the MLE and sample Fr\'echet mean even though the parameter space orbit has to be estimated in the adaptive MRE.  Agreeing with intuition, the adaptive MRE does not perform as well in the case where $n = 5$, $\kappa_i = 2$, and $ \lambda_{ij}  =3$ since in this setting the MLE orbit estimate will not be accurate as the sample size is too small. Moreover, the effect of the choice of orbit on the performance of the adaptive MRE may be larger in this case since $\lambda_{ij}$ is relatively large.

\section{Discussion}
In this work we have shown how estimating the Fr\'echet mean of a distribution can be viewed as an equivariant estimation problem under the isometry group action. Expressions for the minimum risk equivariant estimator are derived in terms of the right Haar measure over the isometry group of the sample space manifold. The MRE can be found explicitly for exponential family models related to the von Mises-Fisher distribution where the density is a decreasing function of the Riemannian distance. These models are parameterized in terms of their Fr\'echet means and a concentration parameter. Developing additional exponential family models on Riemannian manifolds that allow for tractable inferences to be made on the Fr\'echet mean is a future area of interest. It is also of interest to consider models that go beyond random sampling and incorporate dependence between manifold-valued random objects.

When it is not possible to determine the MRE explicitly it can be computed by MCMC over the isometry group, which typically is $\text{O}(k)$ or some product thereof. Monte Carlo methods for simulating manifold-valued random variables is an active area of research \cite{jauch2020random,girolami2011riemann}. The existence of the MRE is guaranteed only when the isometry group acts transitively on the parameter space and the posterior distribution under the right Haar prior has a unique Fr\'echet mean set. We have proposed an adaptive MRE procedure for cases where the action is non-transitive. Non-transitivity is anticipated in settings where the isometry group is ``small", such as the case where the sample space is not homogeneous. For example, stratified spaces, like the BHV space of trees \cite{billera2001geometry}, are not homogeneous but still can have non-trivial isometry groups that can incorporated into an adaptive MRE procedure.

Simulations on the space of positive definite matrices and the torus provide some evidence that the adaptive MRE performs well relative to the estimator it is derived from. It is an open question as to whether the adaptive MRE will always asymptotically dominate any consistent, equivariant estimator from which it is derived.

\section{Appendix}
\subsection{Proofs}
\begin{customlemma}{1}
\label{vonMisLemma}
The Fr\'echet mean under the intrinsic, angular metric of $\mathbb{S}^n$ of the von Mises-Fisher distribution $P_{(\mu,\kappa)}$ is $\mu$ for all $\kappa > 0$. 
\end{customlemma}
\begin{proof}
For $X \sim P_{(\mu,\kappa)}$ we show that $\mathbb{S}^n$, $E(d(X,p)^2) > E(d(X,\mu)^2)$ for any $p \neq \mu$.  As the von Mises-Fisher distributions are invariant, by two-point homogeneity of the sphere it can be assumed without loss of generality that $\mu = (\cos(\theta),\sin(\theta),0,\ldots,0)$ and $p = (\cos(\theta),-\sin(\theta),0,\ldots,0)$ with $\sin(\theta) > 0$. Define $R$ to be reflection about the second coordinate, $R(x_1,x_2,\ldots,x_{n+1}) = (x_1,-x_2,\ldots,x_{n+1})$. Decomposing $\mathbb{S}^n$ into the sets $H_1 \coloneqq \{x \in \mathbb{S}^n: x_2 > 0\}$ and $H_2 = H_1^c$ it is found that
\begin{align*}
E(I_{H_2}\big(d(X,p)^2 - d(X,\mu)^2)\big) & = \int_{H_2} \big(d(x,p)^2 - d(x,\mu)^2\big)\exp(\kappa x^\intercal\mu) \text{vol}(dx)
\\
& =  \int_{H_2} \big(d(Rx,Rp)^2 - d(Rx,R\mu)^2\big)\exp(\kappa (Rx)^\intercal R\mu) \text{vol}(dx) 
\\
& =  \int_{H_1} (d(x,\mu)^2 - d(x,p)^2)\exp(\kappa x^\intercal R\mu) \text{vol}(dx), 
\end{align*}
where a change of variables from $Rx$ to $x$ is used in the last line. Pairing this expectation with the case when $X \in H_1$ gives
\begin{align}
\label{pfofVonmis}
    E\big(d(X,p)^2 - & d(X,\mu)^2 \big) =  \int_{H_1} (d(x,p)^2 - d(x,\mu)^2)\big(\exp(\kappa x^\intercal \mu) - \exp(\kappa x^\intercal R\mu)\big) \text{vol}(dx).
\end{align}
This completes the proof because the integrand above is positive over $H_1$,  since $x^\intercal \mu > x^\intercal R\mu$ whenever $x_2 > 0$. 
\end{proof}

The ideas used in the proof above hold more generally for any two-point homogeneous space.

\begin{customthm}{2}
\label{generalFrecmnthm}
If $\mathcal{X}$ is a two-point homogeneous Riemannian manifold and $X$ is a random object taking values in $\mathcal{X}$ with the density $f\big(d(x,\mu)\big)h(x)$ with respect to the Riemannian volume measure, then $EX = \mu$ if $f$ is a decreasing function. 
\end{customthm}
\begin{proof}
The same proof as in Lemma \ref{vonMisLemma} applies with $H_1 = \{x \in \mathcal{X}: d(x,\mu) < d(x,p)\}$ and $H_2 = H_1^c$. The analogue of the reflection $R$ is an isometry $R$ with $R(p) = \mu, \; R(\mu) = p$. Notice that if $d(x,p) \leq d(x,\mu)$ then $d(Rx,\mu) \leq d(Rx,p)$, implying that $RH_2 = H_1 \bigcup \{x:d(x,p) = d(x,\mu)\}$. As $f(d(x,\mu))  > f(d(x,p))$ in $H_1$, the integral analogous to \eqref{pfofVonmis} is positive as needed. 
\end{proof}

\begin{customlemma}{3}
$G_{\theta_0}$ is compact and
$G$ is closed in $\text{Iso}(\mathcal{X}^n)$.
\label{grouptoplemma}
\end{customlemma}
\begin{proof}
The topology of $G$ is the compact-open topology which means that $g_n \rightarrow g$ if and only if $g_n\vert_K \rightarrow g\vert_K$ uniformly for any any compact $K \subset \mathcal{X}$. 
By a result in \cite{lee2013smooth}, since $G$ is a Lie group acting properly on $\mathcal{X}$, $G_{EP_{\theta_0}}$ is compact. As $G_{\theta_0} \subset G_{EP_{\theta_0}}$ it suffices to show that $G_{\theta_0}$ is closed in $G$. Suppose that $g_n \rightarrow g$ with $g_n \in G_{\theta_0}$ so $P_{g_n\theta_0} = P_{\theta_0}$ for all $n$. As $G$ is a Lie group under the compact-open topology $g_n \rightarrow g$ implies $g_n^{-1} \rightarrow g^{-1}$. For any set $S \subset \mathcal{X}$ let  $S_{\epsilon} = \{x: d(x,S) < \epsilon\}$ be its $\epsilon$-enlargement. By the uniform convergence on compacts, for any $\epsilon$ there exists an $m$ such that $n \geq m$ implies $g_n^{-1}K \subset (g^{-1}K)_{\epsilon}$ and $g^{-1}K \subset (g_n^{-1}K)_{\epsilon}$. This gives the inequalities
\begin{align}
    P_{\theta_0}(K) = P_{\theta_0}(g_n^{-1}K) \leq P_{\theta_0}((g^{-1}K)_{\epsilon}), \label{ineq1}
    \\
 P_{\theta_0}(g^{-1}K) \leq P_{\theta_0}((g_n^{-1}K)_\epsilon) = P_{\theta_0}(g_n^{-1}K_\epsilon)= P_{\theta_0}(K_\epsilon), \label{ineq2}
\end{align}
where the first equality in \eqref{ineq2} holds because each $g_n^{-1}$ is an isometry.
Taking $\epsilon \rightarrow 0$ in \eqref{ineq1},\eqref{ineq2} shows that $P_{\theta}(K) = P_{g\theta_0}(K)$ which implies that $P_{\theta_0} = P_{g\theta_0}$ because $P_{\theta_0}$ and $P_{g\theta_0}$ are inner regular.  

To prove the second statement, let $\Tilde{g}_n \coloneqq g_n \times \cdots \times g_n \rightarrow \Tilde{g}$ in the compact-open topology of $ \text{Iso}(\mathcal{X}^n)$. Then $(g_nx,\ldots,g_nx) \rightarrow \Tilde{g}(x,\ldots,x)$ so $\Tilde{g}(x,\ldots,x) = (y,\ldots,y)$ for some $y \in \mathcal{X}$. Define $g:\mathcal{X} \rightarrow \mathcal{X}$ by $g(x) = y$. It can be checked that $g \in \text{Iso}(\mathcal{X})$ and $g \times \cdots \times g = \Tilde{g}$ so $\Tilde{g} \in G$ as needed.
\end{proof}

\bibliography{biblio}

\end{document}